# Topological invariants of 3-manifolds with boundary

Luca Di Beo

March 2021

## 1 Abstract

This paper presents, with explanatory details, the handle decompositions, fundamental groups and homology groups of 3-manifolds, including some knot complements. Hence, along this paper, when the word *manifold* appears it is implicit that its dimension is 3, except when explicitly generalized for *n* dimensions, $n \in \mathbb{N}$. The results were obtained for: 3-torus ($T^3 = S^1 \times S^1 \times S^1$), projective space $P^3$, trefoil ($3_1$), figure-eight ($4_1$), cinquefoil ($5_1$) and three-twist ($5_2$).

## 2 Introduction: Preliminary definitions

The fundamental group and homology group theories are general enough to comprise manifolds in any dimension. However, these topological invariants, and others, are by far most practical and widely spread used in the context of 3-manifolds. This is so since every person is surrounded by them in the day-to-day life – as the physical space presents itself in nature in 3 dimensions.
With this motivation, the technique known as Heegaard splitting, for computing handle decompositions, will be used here to simplify the calculation of topological invariants of *n*-manifolds (with emphasis on *n* = 3).

### 2.1 Manifolds with boundary [3] [4]

A Hausdorff topological space *M* with a countable base is called a *n-dimensional manifold with boundary* if each of its points has a neighbourhood homeomorphic either to the space $\mathbf{R}^n = \{(x_1, ..., x_n) \mid x_i \in \mathbf{R}\}$ (i.e., a ball neighbourhood) or to the closed subspace $\mathbf{\bar{R}}^n = \{(x_1, ..., x_n) \in \mathbf{R}^n \mid x_n \geq 0\}$: (i.e., a half ball neighbourhood).[1] The union of the points of the manifold *M*, which do not admit ball neighbourhoods (but perhaps only half ball neighbourhoods), is called the boundary of this manifold and is denoted by $\partial M$. The figure below (one of the famous drawings by Professor Anatoly Fomenko, of the Moscow State University) shows a few examples of manifolds with boundary:

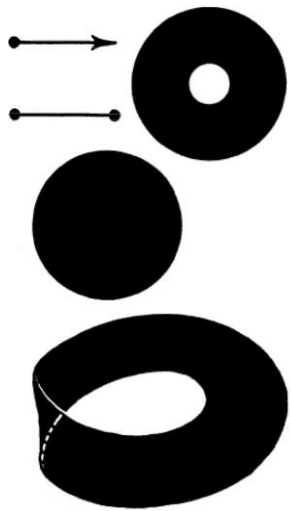

---

[1]The conditions *countable base*, i.e. the space can be covered by a countable number of ball neighbourhoods, and to be *Hausdorff*, i.e. for any two points of the the space there must exist disjoint neighbourhoods, can be dropped with the consequence of pathological cases.

### 2.1.1 Local coordinate systems about boundary points [3]

Let $M$ be a $n$-manifold with boundary and $p \in \partial M$. Then there exist a local coordinate system $(x_1, ..., x_n)$, with $x_n \geq 0$, about $p$. Hence, $M$ corresponds to $x_n \geq 0$ and $\partial M$ to $x_n = 0$.

### 2.1.2 Construction of manifolds [3]

Consider $M$ as an $n$-dimensional manifold (without boundary) and let $f : M \to R$ be a smooth function ($C^\infty$) on $M$. Assume that 0 is not a critical value[2]. Define a subset of $M$ as $M_{f \geq 0} = \{p \in M \mid f(p) \geq 0\}$ – thus $M_{f \geq 0}$ is a $n$-dimensional manifold with boundary. Its boundary is $\partial M_{f \geq 0} = M_{f=0} = \{p \in M \mid f(p) = 0\}$. In fact, the disk $D^n$ (see the first example in subsection 2.1.2) can be defined by taking $R^n$ as $M$ and defining $f : R^n \to R$ by $f = 1 - (x_1^2 + ... + x_n^2)$.

### 2.1.3 Examples of manifold with boundary [3]

**The $n$-disk:**

The Euclidean space of dimension $n$ can be represented by $R^n$. The *n-dimensional (unit) disk* (*n*-disk) is defined as $D^n = \{(x_1, ..., x_n) \in R^n \mid x_1^2 + ... + x_n^2 \leq 1\}$. The $n$-disk is an example of $n$-dimensional manifold with boundary. Its boundary is the $(n-1)$-dimensional unit sphere, defined as $S^{n-1} = \{(x_1, ..., x_n) \in R^n \mid x_1^2 + ... + x_n^2 = 1\}$. Alternatively, the boundary of $D^n$ is denoted as $\partial D^n = S^{n-1}$. On account of the focus of this paper, which is to describe manifolds in 3 dimensions, it is depicted below the 3-disk $D^3$ (a 3-manifold with 2-dimensional boundary $S^2$, i.e., a sphere) – which will be an useful illustration to keep in mind when studying 0- and 3-handles in 3-manifolds:

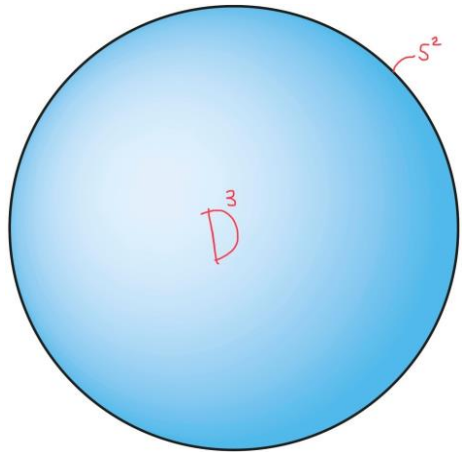

**The $n$-dimensional upper half space:**

The set $R_+^n$ is a $n$-dimensional manifold with boundary. Its boundary is defined as $R^{n-1} = \{(x_1, ..., x_{n-1}) \in R^{n-1}\}$ – denoted as $\partial R^n_+ = R^{n-1}$. Once again, the 3-dimensional case can be easily visualized:

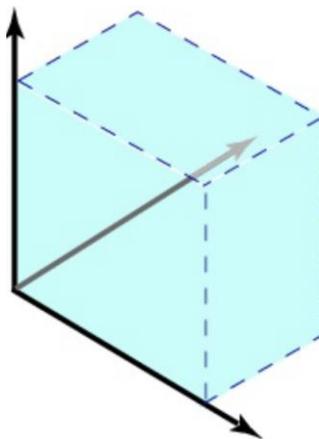

---

[2]A point $p$ is called a *critical point* of $f$ if $\frac{\partial f}{\partial x_i}(p) = 0$ for $i = 1, ..., n$. of $f$.

### 2.1.4 Morse functions

Morse functions are important for definitions (and propositions) in the realm of handle decomposition of manifolds, and therefore will be defined: If none of the critical points of a function $f : M \to R$ are degenerate[3], then $f$ is called a *Morse function*.

### 2.1.5 A special case of the Van Kampen's theorem [3]

This special case will be considered for 3-dimensional spaces only. For a general definition of Van Kampen's theorem, consult [5]. Let $h : \partial D^2 = S^1 \to X$ be a continuous mapping with codomain $X$ (3-dimensional space). The space $X \cup_h D^2$ is obtained by attaching $D^2$, by its boundary, to $X$ through the mapping $h$. Setting $x_0$ as a base point of $X$, it is easily seen that the inclusion mapping $i : (X, x_0) \hookrightarrow (X \cup_h D^2, x_0)$ induces the homomorphism $i_* : \pi_1(X, x_0) \hookrightarrow \pi_1(X \cup_h D^2, x_0)$, which is bijective, with kernel $Ker(i_*) \equiv N(\{h\})$, where $N(\{h\})$ is the normal subgroup of $\pi_1(X, x_0)$ generated by $\{h\}$.[4]

Given the surjectivity of $i_*$, the fundamental theorem on homomorphisms [2] implies that:

$$\pi_1(X \cup_h D^2, x_0) \cong \pi_1(X, x_0)/N(h)$$

## 2.2 Handle decomposition of manifolds

### 2.2.1 Definition of Handles [1]

A *k-handle* is defined as the subset of a topological space $M$ that is homeomorphic to $D^k \times D^{n-k}$, with $k \in N$, $k \leq n$ and $n = dim(M)$. There are several important subsets of the $k$-handle:

- The *core disk* – $D^k \times 0$;
- The *attaching sphere* – $S^{k-1} \times 0$;
- The *co-core disk* – $0 \times D^{n-k}$;
- The *belt sphere* – $0 \times S^{n-k-1}$;
- The *attaching region* – $S^{k-1} \times D^{n-k}$;
- The *belt region* – $D^k \times S^{n-k-1}$.

As easily seen, every $k$-handle is also an $(n-k)$-handle. However, a distinction between them will rise during the process of gluing, according to the order of the disks defining the expression $D^k \times D^{n-k}$. It is standard to considered the second disk, i.e. $D^k$, as the one that is glued to the space $M$. This identification produces a deference between gluing a $k$-handle and a $(n-k)$-handle to $M$, as one can easily check by the first figure in section 2.2.2, which is another drawing by Professor Anatoly Fomenko, of the Moscow State University. The process of transforming a $k$-handle into a $(n-k)$-handle is usually expressed as the: "handle decomposition was turned upside-down".

Every handle is homeomorphic to an $n$-disk. Handles in various dimensions, as well as important spheres and disks, are presented below. Especially important for this paper are the 3-dimensional cases.

---

[3]A degenerate point in this context is understood as a point that has no unique tangent passing through it.
[4]$\{h\}$ is the set of loops.

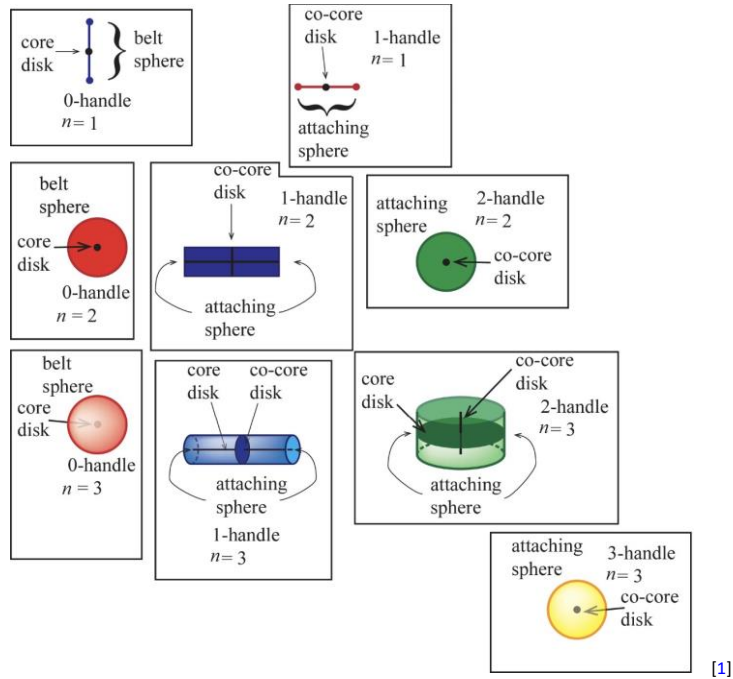

[1]

### 2.2.2 Handlebody [3]

A *handlebody* is a manifold with boundary obtained from $D^n$ by attaching $m$ handles of various indices one after another, and is denoted as $H(D^n; \varphi_1, ..., \varphi_m) := D^n \cup_{\varphi_1} D^{k_1} \times D^{n-k_1} \cup_{\varphi_2} ... \cup_{\varphi_m} D^{k_m} \times D^{n-k_m}$.

These attachings are realized through a $C^\infty$ mapping $\varphi_i : \partial D^{k_i} \times D^{n-k_i} \to \partial D^n$, with $i = 1, ..., m$. The figure below shows the attaching of *0,1,2*-handles to an arbitrary manifold with boundary:[5] [4]

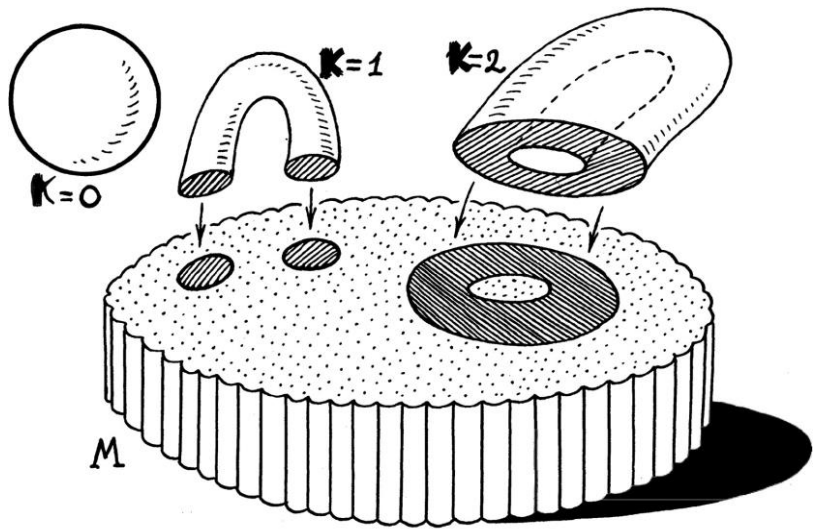

### 2.2.3 3-dimensional handlebody and Heegaard splitting [3]

Let $M$ be a connected, orientable, closed 3-manifold. Let $f : M \to \mathbb{R}$ be a Morse function such that the handlebody decomposition associated with $f$ is $M := h^0 \cup (h^1_1 \cup ... \cup h^1_{k_1}) \cup (h^2_1 \cup ... \cup h^2_{k_2}) \cup h^3$, where $h_j$ represents a *j*-handle ($j \in \mathbb{N}$). In a 3-dimensional handlebody, as depicted below, $H_g$ is called *handlebody of genus g* – in the case of this figure $g = 3$. Its boundary $\partial H_3$ is an orientable closed surface of genus 3.

---

[5]The attaching of the 0-handle to the manifold is simply a disjoint union.

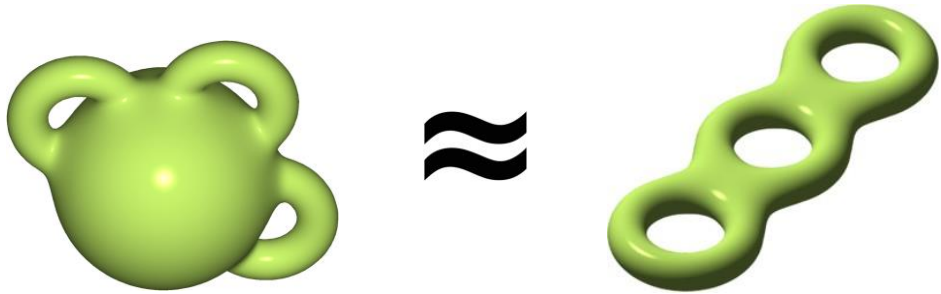

Let $H$ be a subhandlebody of the handle decomposition $H_3 = h^0 \cup (h^1_1 \cup h^1_2 \cup h^1_3) \cup h^3$ above and $H^* := H_3 - int(H)$ – then $H^*$ is the union of a 0-handle and 1-handles. Thus, $H_3$ can be described as the union $H_3 = H \cup H^*$, which is obtained by attaching the handlebodies $H$ and $H^*$ along their boundaries by the mapping $\varphi : \partial H^* \to \partial H$.

This description of $H_3$ is called a *Heegaard splitting* (*of genus 3* in this case). The belt spheres of 1-handles ($0 \times S^1$) contained in the arbitrary handlebody $H$ of genus $k$ form $k$ copies of circles in the boundary $\partial H$. This circle is called a *meridian*. They are denoted by $m_1, ..., m_k$. Similarly, meridians $m^*_1, ..., m^*_k$ are defined for $H^*$.

In this Heegaard splitting of $H_3$, the images in $\partial H$ of the meridians of $H^*$, under the mapping $\varphi$ are denoted as $\mu_1 = \varphi(m^*_1), ..., \mu_k = \varphi(m^*_k)$.[6]

Formally, the handlebody $H$ of genus $k < 3$, together with the closed curves $\mu_1, ..., \mu_k$ on the boundary $\partial H$ obtained above, namely $(H; \mu_1, ..., \mu_k)$, is called a *Heegaard diagram* of $H_3$.

Since handle decompositions of $H_3$ are not unique, Heegaard diagrams are not as well – remember that $H$ was picked arbitrarily.

### 2.2.4 General example of handlebody

Before moving on to calculate the fundamental and homology groups of knot complements, it will be presented an example of these calculations applied to a handlebody $T^3$ of genus 3 – namely, the 3-torus:

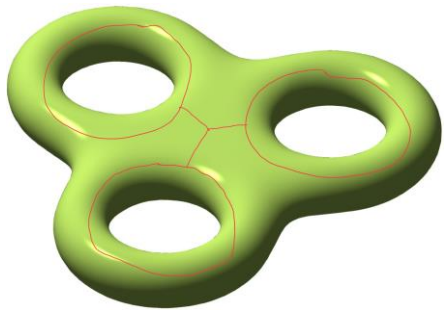

The three connected red circles drawn on the figure above are homeomorphic to $S^1 \times S^1 \times S^1$. Naming each one of these circles as $a$, $b$ or $c$, it is possible to calculate the fundamental group of $T^3$ as a presentation in terms of generators and relations:

$$\pi_1(T^3) = \pi_1(S^1 \times S^1 \times S^1) = F^3 = <a, b, c \mid >$$

As seen, there are no relations. The result above can be obtained through the traditional Van Kampen's theorem [5] (not the special case shown previously here).

The homology groups are as follow:

$$H_1(T^3) = \mathbb{Z} \oplus \mathbb{Z} \oplus \mathbb{Z}$$
$$H_0(T^3) = \mathbb{Z}$$
$$H_2(T^3) = 0$$

The fact that the first homology group is homeomorphic to $\mathbb{Z} \oplus \mathbb{Z} \oplus \mathbb{Z}$ reflects the 3 genus topological invariance of this space. As it is clearly seen in the image above, $T^3$ is a neighborhood of $S^1 \times S^1 \times S^1$.

---

[6]This notation shows that the meridians of the handlebody are analysed from the point of view of the $H^*$ rather than of $H$.

*Remark*: The homology group (of a general space $X$) is calculated through the formula

$$H_n(X) \equiv \frac{Ker(\partial_n)}{Im(\partial_{n+1})},$$

where $\partial_n$ are the homomorphisms involved.

## 3  Handle decomposition of projective space $P^3$ [3]

It will be presented a Heegaard diagram of the projective space $P^3$. This space has a handle decomposition of the form $P^3 = h^0 \cup h^1 \cup h^3 \cup h^3$. Suppose that $H \subset P^3$ is a handlebody which is homeomorphic to the *solid torus* ($S^1 \times D^2$). Let $\tau$ denote a fixed point on its core-disk $S^1$, then a meridian will be represented as $m := \{\tau\} \times \partial D^2$. A *longitude* on $\partial H$ is defined as $l := S^1 \times \{\sigma\}$, where $\sigma$ is a fixed point on $\partial D^2$, and intersects $m$ in a single point transversely.

For $x_* := (\tau, \sigma)$ a fixed point on $H$, $\pi_1(H, x_*) \cong \mathbb{Z}$. The cell decomposition of $P^2$ and $P^3$, where $e^j$ is a $j$-dimensional cell, $\forall j \in \mathbb{N}$, are:

$$P^2 = e^0 \cup e^1 \cup e^2$$
$$P^3 = e^0 \cup e^1 \cup e^2 \cup e^3$$

Using the special case of Van Kampen's theorem

$$\pi_1(P^3, x_*) \cong \pi_1(P^2, x_*) \cong \mathbb{Z}/2\mathbb{Z}$$

Therefore, a simple closed curve $\mu$ that provides a Heegaard diagram of $P^3$ has to intersect $m$ twice transversely in the same direction, as depicted in the figure below:

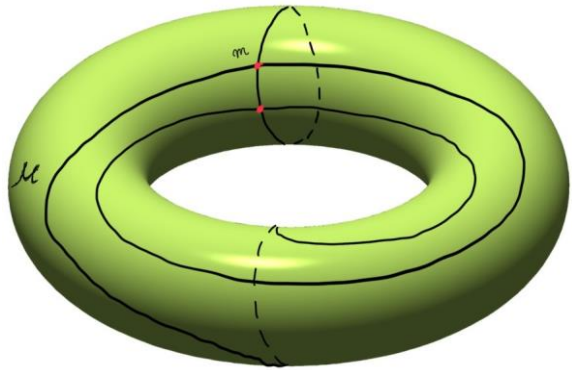

## 4  Knot and knot complement

### 4.1  Construction of knot complement

Consider a mapping $K : S^1 \hookrightarrow S^3$, which is a smooth embedding of a circle into the 3-sphere – $S^3 = \{(x_0, x_1, x_2, x_3) \in \mathbb{R}^4 \mid x_0^2 + x_1^2 + x_2^2 + x_3^2 = 1\}$. The image of this mapping is simply called $K$, and it is referred to as the knot $K$ throughout all the paper.
Consider now a tubular neighborhood of $K$ as a smooth embedding of the solid torus $N : S^1 \times D^2 \hookrightarrow S^3$. See the illustration below:

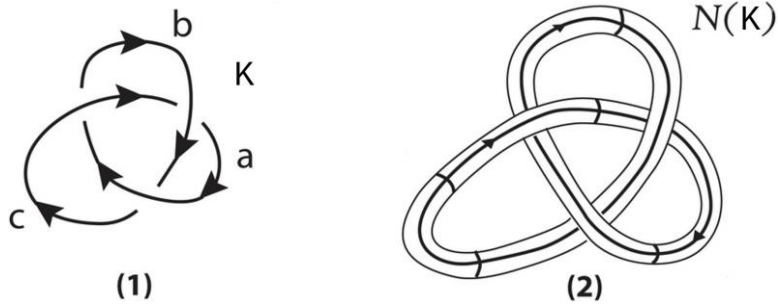

(1)            (2)

Therefore, the knot is embedded as the core $0 \times S^1$ of the solid torus.

The knot complement (also called knot exterior) is the space defined as $E := E(K) := S^3 \setminus int(N)$.

## 4.2 Trefoil: $3_1$

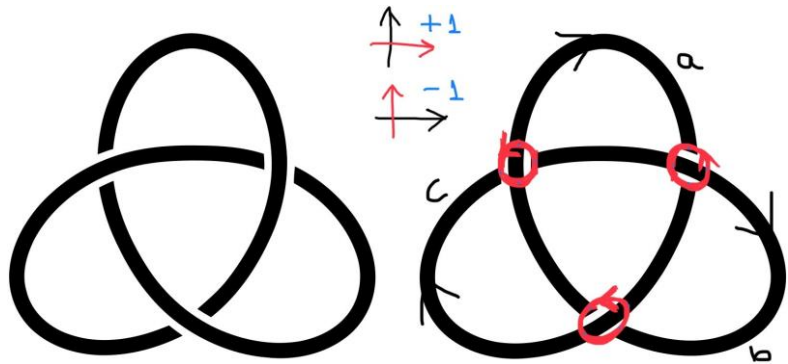

In the figure above (left) the $K$ knot $3_1$ was projected in the plane with the crosses, or intersections, of two arcs of the knot represented by breaking the arc that is closer to the plane of projection than an observer is – this is a standard representation, and therefore will be used for every knot in this paper. Each arc was named *a, b* or *c*. At the right side of the figure above, 2-handles were attached by their boundaries (red circles) to each intersection[7]. Arrows were placed in each arc, and in each boundary of the 2-handles, in such a way to provide consistent orientations. Moreover, in the top middle of the figure above, a system of orientation was defined (+1 or -1) – *which will be the same used in all the following sections*.

With all these tools well established, it is possible to calculate the fundamental group and, subsequently, the homology group of the $3_1$ complement:

$$\pi_1(S^3 \setminus K) = < a, b, c \mid a^{-1}cab^{-1} \equiv e, abc^{-1}b^{-1} \equiv e, ca^{-1}c^{-1}b \equiv e >$$

Notice that these three relations imply, respectively, that (I) $a \equiv b$, (II) $a \equiv c$ and (III) $a \equiv b$, hence (I) and (II) $\Rightarrow$ (III), i.e., (III) can be omitted since it does not provide any new information – *as it will be omitted in all the other examples of this paper from now on*. Thus, the fundamental group (followed by the zero-homology group) is rather written as:

$$\pi_1(S^3 \setminus K) = < a, b, c \mid a^{-1}cab^{-1} \equiv e, abc^{-1}b^{-1} \equiv e >$$

$$H_0(S^3 \setminus K) = \mathbb{Z}$$

This is so since $K$ is a connected manifold, i.e., isomorphic to the circle.

$$H_1(S^3 \setminus K) = \frac{\pi_1}{[\pi_1, \pi_1]} = \frac{< a, b, c \mid a = b = c >}{0} = < a > \approx \mathbb{Z}$$

It turns out that the first homology group of every knot complement will be $\mathbb{Z}$.

---

[7]The same subsequent results could be achieved by attaching 1-handles to each intersection as well.

$$0 \xrightarrow{\partial_3} C_2 \xrightarrow{\partial_2} C_1 \xrightarrow{\partial_1} C_0 \xrightarrow{\partial_0} 0 \qquad \partial_0 = \partial_1 = \partial_3 = 0$$

$$C_2 = <a^{-1}cab^{-1}, abc^{-1}b^{-1}>$$
$$C_1 = <a, b, c>$$
$$C_0 = <*> \qquad \text{(generated by a point)}[8]$$

$$\partial_2(a^{-1}cab^{-1}) = c - b \qquad \partial_2(abc^{-1}b^{-1}) = a - c$$
$$\text{Ker}(\partial_2) = 0 \implies H_2(S^3 \setminus K) = \frac{\text{Ker}(\partial_2)}{\text{Im}(\partial_3)} = 0$$

The second homology group of any knot complement is always zero.
Hence, the calculation of the homology groups of the subsequent knots will be the same, and their fundamental group will be calculated in a similar way as well.

### 4.3 Figure-eight: $4_1$

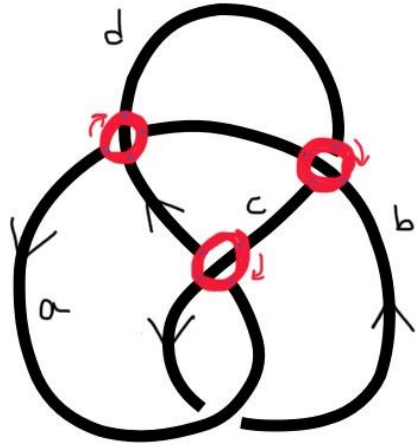

The $K$ knot depicted here is the $4_1$.
As before the fundamental and homology groups are as follow:

$$\pi_1(S^3 \setminus K) = <a, b, c, d \mid db^{-1}d^{-1}a \equiv e, d^{-1}b^{-1}cb \equiv e, c^{-1}a^{-1}cd \equiv e>$$

$$C_2 = <db^{-1}d^{-1}a, d^{-1}b^{-1}cb, c^{-1}a^{-1}cd>$$
$$C_1 = <a, b, c, d>$$
$$C_0 = <*>$$

$$\partial_2(db^{-1}d^{-1}a) = a - b \qquad \partial_2(d^{-1}b^{-1}cb) = c - d \qquad \partial_2(c^{-1}a^{-1}cd) = d - a$$
$$H_0(S^3 \setminus K) = \mathbb{Z}$$
$$H_1(S^3 \setminus K) = \mathbb{Z}$$
$$H_2(S^3 \setminus K) = 0$$

---
[8]I.e. this handle is the neighborhood of a 0-cell.

## 4.4 Cinquefoil: $5_1$

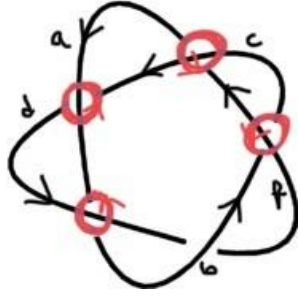

$\pi_1(S^3 \setminus K) = < a, b, c, d, f \mid ad^{-1}a^{-1}c \equiv e, fca^{-1}c^{-1} \equiv e, fc^{-1}f^{-1}b \equiv e, b^{-1}d^{-1}ad \equiv e >$

$C_2 = < ad^{-1}a^{-1}c, fca^{-1}c^{-1}, fc^{-1}f^{-1}b, b^{-1}d^{-1}ad >$
$C_1 = < a, b, c, d, f >$
$C_0 = < * >$

$\partial_2(ad^{-1}a^{-1}c) = c - d \qquad \partial_2(fca^{-1}c^{-1}) = f - a \qquad \partial_2(fc^{-1}f^{-1}b) = b - c \qquad \partial_2(b^{-1}d^{-1}ad) = a - b$

$H_0(S^3 \setminus K) = \mathbb{Z}$
$H_1(S^3 \setminus K) = \mathbb{Z}$
$H_2(S^3 \setminus K) = 0$

## 4.5 Three-twist: $5_2$

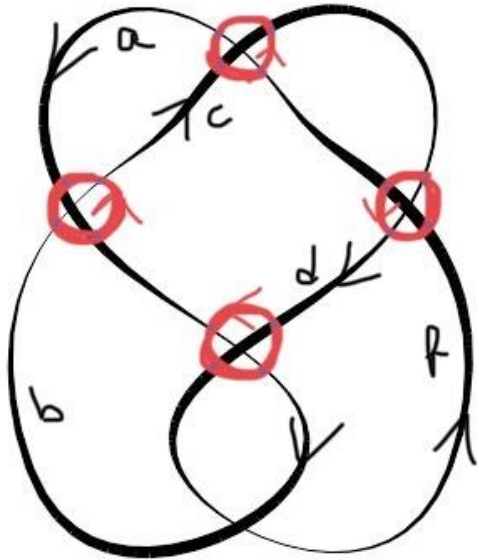

$\pi_1(S^3 \setminus K) = < a, b, c, d, f \mid a^{-1}cfc^{-1} \equiv e, c^{-1}aba^{-1} \equiv e, dad^{-1}b^{-1}b \equiv e, cf^{-1}d^{-1}f \equiv e >$

$C_2 = < a^{-1}cfc^{-1}, c^{-1}aba^{-1}, dad^{-1}b^{-1}b, cf^{-1}d^{-1}f >$
$C_1 = < a, b, c, d, f >$
$C_0 = < * >$

$\partial_2(a^{-1}cfc^{-1}) = f - a \qquad \partial_2(c^{-1}aba^{-1}) = b - c \qquad \partial_2(fc^{-1}f^{-1}b) = a - b \qquad \partial_2(dad^{-1}b^{-1}b) = c - d$

$H_0(S^3 \setminus K) = \mathbb{Z}$
$H_1(S^3 \setminus K) = \mathbb{Z}$
$H_2(S^3 \setminus K) = 0$

# 5 Conclusion

Heegaard diagrams are useful for simplifying the computation of topological invariants – in particular the fundamental and homology groups – of 3-manifolds with boundary. Studying their essential features is extremely important for applied mathematics, physics, chemistry and biology, since every person is surrounded by several examples of 3-manifolds in their day-to-day experiences. Therefore, 3-manifolds and their invariants have direct influences in processes concerning the fields cited above.